\documentclass{amsart}

\usepackage{amsfonts, amssymb}
\usepackage{amsmath}
\usepackage{amsthm}

\numberwithin{equation}{section}

\newtheorem{theorem}{Theorem}[section]

\newtheorem{lemma}[theorem]{Lemma}
\newtheorem{example}[theorem]{Example}

\newtheorem{remark}[]{Remark}

\begin{document}

\title[solutions of  nonlinear
differential-difference equations]{Some results on transcendental entire solutions of certain nonlinear
differential-difference equations}

\author{Nan Li}
\address{School of Mathematics, Qilu Normal University, Jinan, Shandong, 250200, P.R. China}
\email{nanli32787310@163.com}

\author{Jiachuan Geng}
\address{School of Public Administration and Policy, Shandong university of Finance and Economics,
  Jinan, Shandong, 250014, P.R.China}
\email{cosart@126.com}

\author{Lianzhong Yang}

\address{School of Mathematics, Shandong University,
  Jinan, Shandong, 250100, P.R.China}
\email{lzyang@sdu.edu.cn}


\thanks{ This work was supported by  NNSF of China (No.11801215), and the NSF of Shandong Province, P. R. China (No.ZR2016AQ20 \& No. ZR2018MA021).}

\subjclass[2010]{Primary 30D35; Secondary 39B32}

\date{\today}

\commby{}

\begin{abstract}
In this paper, we study the  transcendental entire solutions
  for the nonlinear differential-difference equations of the forms:
  \begin{eqnarray*}
    f^{2}(z)+\widetilde{\omega} f(z)f'(z)+q(z)e^{Q(z)}f(z+c)=u(z)e^{v(z)},
  \end{eqnarray*}
  and
  \begin{eqnarray*}
       f^{n}(z)+\omega f^{n-1}(z)f'(z)+q(z)e^{Q(z)}f(z+c)=p_{1}e^{\lambda_{1} z}+p_{2}e^{\lambda_{2} z}, \quad n\geq 3,
      \end{eqnarray*}
 where $\omega$ is a constant, $\widetilde{\omega}, c, \lambda_{1}, \lambda_{2}, p_{1}, p_{2}$ are non-zero constants,  $q, Q, u, v$ are polynomials such that $Q,v$ are not constants and $q,u\not\equiv0$. Our results are improvements and complements of some previous results.
\end{abstract}

\maketitle

\section{Introduction}

Let $f(z)$ be a transcendental meromorphic function in the complex plane
$\mathbb{C}$.  We assume that the reader is familiar with the standard notations and main results in
Nevanlinna  theory (see \cite{Hayman},\cite{Laine},\cite{Yi1}).
Throughout this paper, the term $S(r,f)$ always has the property that
$S(r,f)=o(T(r,f))$ as $r\to \infty$, possibly outside a set $E$ (which is not necessarily the same at each occurrence) of finite
logarithmic measure. A meromorphic function $a(z)$ is said to be a small function
with respect to $f(z)$ if and only if $T(r,a)=S(r,f)$.

In the past few decades, many scholars have investigated existence, order and value distribution  of solutions of complex differential or difference equations, see \cite{chen2011, Gundersen, Laine, Long} etc.

In 1964, Hayman \cite{Hayman} considered the following non-linear differential equation
\begin{eqnarray}\label{thmABaddAeq1}
  f^{n}(z)+Q_{d}(f(z))=g(z),
\end{eqnarray}
where $Q_{d}(f)$ is a differential polynomial in $f$ with
degree $d$ and obtained the following
result.

\begin{theorem}\label{thmABaddA}
Suppose that $f(z)$ is a nonconstant meromorphic function, $d\leq n-1$, and
$f,g$ satisfy $N(r,f)+N(r,1/g)=S(r,f)$ in \eqref{thmABaddAeq1}. Then we have
$g(z)=(f(z)+\gamma(z))^{n}$, where $\gamma(z)$ is meromorphic
and a small function of $f(z)$.
\end{theorem}

After that, dozens of  papers , see \cite{LiP2011,Liao2013,LIAO2015,Lu2017} etc., focus on the solutions of the nonlinear differential equations of the forms
\begin{equation*}
  f^{n}(z)+P(f(z)) =h(z), \; \textrm{or} \;  f^{n}(z)f'(z)+P(f(z)) =h(z),
\end{equation*}
where $P(f)$ denotes a differential polynomial in $f$
of degree at most $n-1$, and
$h$ is a given  meromorphic function such as $ue^{v}$(where $u,v$ are polynomials),
$p_{1}e^{\alpha_{1}}+p_{2}e^{\alpha_{2}}$(where $p_{1}, p_{2}, \alpha_{1},\alpha_{2}$ are polymomials), etc.

In 2012, Wen et al. \cite{Wen2012} investigated and classified the finite order entire solutions of
the equation
\begin{eqnarray}\label{introeq1zzz}
 f^{n}(z)+q(z)e^{Q(z)}f(z+c)=P(z),
\end{eqnarray}
where $q,Q,P$ are polynomials, $n\geq 2$ is an integer, and $c\in \mathbb{C}\setminus\{0\}$. Later, Chen \cite{chen2018JCAA} replaced $P(z)$ in \eqref{introeq1zzz} by $p_{1}e^{\lambda z}+p_{2}e^{-\lambda z}$, where $p_{1},p_{2},\lambda$ are non-zero constants, and studied its finite order entire solutions when $n\geq3$.

By observing all the above equations, we can see that there exists only one dominant term $f^{n}$ or $f^{n}f'$ on the  left-hand side. Hence
a natural interesting area of inquiry
is the study of equations which may have two dominated terms on the left-hand side with the same degree.

Motivated by the above equations, Chen, Hu and Wang \cite{Chen-Hu-Wang} investigated the following
non-linear differential-difference equation
\begin{equation}\label{introeq1add2}
  f^{n}(z)+\omega f^{n-1}(z)f'(z)+q(z)e^{Q(z)}f(z+c)=u(z)e^{v(z)},
\end{equation}
where $n$ is a positive integer, $c\neq 0$, $\omega$ are constants, $q,Q,u,v$ are polynomials such that  $Q, v$ are not constants and $q,u \not\equiv0$,  and obtained the following result.

\begin{theorem}\label{thmAB}
Let $n$ be an integer satisfying $n\geq 3$ for $\omega\neq0$ and
$n\geq2$ for $\omega=0$. Suppose that $f$ is a non-vanishing transcendental
entire solution of finite order of \eqref{introeq1add2}. Then
every solution $f$ satisfies one of the following results:
\begin{itemize}
  \item [(1)] $\rho(f)<\deg v=\deg Q$ and $f=Ce^{-z/\omega}$, where
  $C$ is a constant.
  \item [(2)] $\rho(f)=\deg Q\geq \deg v$.
\end{itemize}
\end{theorem}

Then it's natural to ask: what will happen for the solutions of equation
 \eqref{introeq1add2} when $n=2$ and $\omega\neq0$? In this paper, we study this problem and obtain the
 following result, which is a complement of Theorem~\ref{thmAB}.

\begin{theorem}\label{thmABABnew1}
Let $c,\, \widetilde{\omega}\neq 0$ be constants, $q, Q, u, v$ be  polynomials such that $Q,v$ are not constants and $q,u\not\equiv0$.  Suppose that $f$ is a  transcendental entire solution
with finite order of
\begin{eqnarray}\label{aaa111}
f^{2}(z)+\widetilde{\omega} f(z)f'(z)+q(z)e^{Q(z)}f(z+c)=u(z)e^{v(z)},
\end{eqnarray}
satisfying $\lambda(f)<\rho(f)$, then $\deg Q=\deg v$, and
one of the following relations holds:
 \begin{itemize}
   \item [(1)] $\sigma(f)< \deg Q=\deg v$, and $f=Ce^{-z/\widetilde{\omega}}$
   \item [(2)]  $\sigma(f)= \deg Q=\deg v$.
 \end{itemize}
\end{theorem}

Two examples are given below to show that our estimates in Theorem~\ref{thmABABnew1} are sharp.

\begin{example}\label{example1}
$f_{0}(z)=2e^{-z}$ is a transcendental entire solution of the following
differential-difference equation
\begin{eqnarray*}
  f^{2}+ff'+ze^{z^{2}+z+1}f(z+1)=2ze^{z^{2}}.
\end{eqnarray*}
Here $\widetilde{\omega}=1\neq0$, $Q=z^{2}+z+1$, $v=z^{2}$ and $0=\lambda(f_{0})<\sigma(f_{0})=1$. Then we have $\sigma(f_{0})=1<2=\deg Q=\deg v$, and $f_{0}=Ce^{-z/\widetilde{\omega}}$, where $C=2$. This shows that Conclusion (1) of
Theorem~\ref{thmABABnew1} occurs.
\end{example}

\begin{example}\label{example1}
$f_{1}(z)=e^{z^{2}}$ is a transcendental entire solution of the following
differential-difference equation
\begin{eqnarray*}
  f^{2}+ff'+e^{z^{2}-2z-1}f(z+1)=2(z+1)e^{2z^{2}}.
\end{eqnarray*}
Here $\widetilde{\omega}=1\neq0$, $Q=z^{2}-2z-1$, $v=2z^{2}$ and $0=\lambda(f_{1})<\sigma(f_{1})=2$.
Then we have $\sigma(f_{1})=2=\deg Q=\deg v$. This illustrates that Conclusion (2) of
Theorem~\ref{thmABABnew1} also exits.
\end{example}

It is also interesting to investigate the entire solutions with
finite order of the following differential-difference equation
\begin{eqnarray}\label{thmABABnew2eq1}
  f^{n}(z)+\omega f^{n-1}(z)f'(z)+q(z)e^{Q(z)}f(z+c)=p_{1}e^{\lambda_{1} z}+p_{2}e^{\lambda_{2} z},
\end{eqnarray}
where $n$ is a positive integer, $\omega$ is a constant and $c, \lambda_{1}, \lambda_{2}, p_{1}, p_{2}$ are non-zero constants,   $q,Q$ are polynomials such that  $Q$ is not a constant and $q\not\equiv0$.

In 2020, Chen, Hu and Wang \cite{Chen-Hu-Wang} studied the special case $\lambda_{2}=-\lambda_{1}$ for the above problem and obtained the following
result.
\begin{theorem}\label{thmBA}
If $f$ is a transcendental entire solution with finite order of
\begin{eqnarray}\label{letterthmBAeq1}
 f^{n}(z)+\omega f^{n-1}(z)f'(z)+q(z)e^{Q(z)}f(z+c)=p_{1}e^{\lambda z}+p_{2}e^{-\lambda z},
\end{eqnarray}
then
the following conclusions hold:
\begin{itemize}
   \item [(i)] If $n\geq4$ for $\omega\neq0$ and $n\geq3$ for $\omega=0$, then every solution $f$ satisfies $\rho(f)= \deg Q=1$.
   \item [(ii)]  If $n\geq1$ and $f$ is a solution of \eqref{letterthmBAeq1} which belongs to $\Gamma_{0}$, then
       \begin{eqnarray*}
         f(z)=e^{\lambda z/n+B},\; Q(z)=-\frac{n+1}{n}\lambda z+b
       \end{eqnarray*}
       or
       \begin{eqnarray*}
         f(z)=e^{-\lambda z/n+B},\; Q(z)=\frac{n+1}{n}\lambda z+b
       \end{eqnarray*}
       where $b,B\in \mathbb{C}$, and $\Gamma_{0}=\{e^{\alpha(z)}:\alpha(z)$ is a non-constant polynomial$\}$.
 \end{itemize}

\end{theorem}

For the case $\lambda_{2}=\lambda_{1}$, equation \eqref{thmABABnew2eq1}
can be reduced to $f^{n}(z)+\omega f^{n-1}(z)f'(z)+q(z)e^{Q(z)}f(z+c)=(p_{1}+p_{2})e^{\lambda_{1} z}$, so
we can get the order of entire solutions by using Theorems~\ref{thmAB} and ~\ref{thmABABnew1}.

Then it is natural to ask : what will happen for the entire solutions of
equation \eqref{thmABABnew2eq1} when $\lambda_{2}\neq \pm\lambda_{1}$ ?
In this paper, we study  this problem and obtain the following result.

\begin{theorem}\label{thmABABnew2}
If $f$ is a transcendental entire solution with finite order of
\eqref{thmABABnew2eq1},  then
the following conclusions hold:
\begin{itemize}
   \item [(1)] If $n\geq4$ for $\omega\neq0$ and $n\geq 3$ for $\omega=0$, then every solution $f$ satisfies $\sigma(f)= \deg Q=1$.
   \item [(2)]  If $n\geq1$ and $f$ is a solution of \eqref{thmABABnew2eq1} with
   $\lambda(f)<\sigma(f)$, then
       \begin{eqnarray*}
        f(z)=\left(\frac{p_{2}n}{n+\omega\lambda_{2}} \right)^{\frac{1}{n}}e^{\frac{\lambda_{2}z}{n}},\; Q(z)=\left(\lambda_{1}-\frac{\lambda_{2}}{n}\right)z+b_{1},
       \end{eqnarray*}
 or
       \begin{eqnarray*}
         f(z)=\left(\frac{p_{1}n}{n+\omega\lambda_{1}} \right)^{\frac{1}{n}}e^{\frac{\lambda_{1}z}{n}},\; Q(z)=\left(\lambda_{2}-\frac{\lambda_{1}}{n}\right)z+b_{2},
 \end{eqnarray*}
       where $b_{1}, b_{2}\in \mathbb{C}$ satisfy $p_{1}=q\left(\frac{p_{2}n}{n+\omega\lambda_{2}} \right)^{\frac{1}{n}}e^{\frac{\lambda_{2}c}{n}+b_{1}}$
       and $p_{2}=q\left(\frac{p_{1}n}{n+\omega\lambda_{1}} \right)^{\frac{1}{n}}e^{\frac{\lambda_{1}c}{n}+b_{2}}$, respectively.
 \end{itemize}

\begin{remark}
Obviously, let $\lambda_{2}=-\lambda_{1}$, we can see that Theorem~\ref{thmABABnew2}(1) is an improvement of Theorem~\ref{thmBA} (i).
Since any function $f$ which belongs to $\Gamma_{0}$ all satisfy the condition that
$\lambda(f)<\sigma(f)$, thus Theorem~\ref{thmABABnew2}(2) is also an improvement of Theorem~\ref{thmBA} (ii).

\end{remark}

\end{theorem}

By observing Conclusion (1) of Theorem~\ref{thmABABnew2}, it is natural to
ask what will happen for the entire solutions of
equation \eqref{thmABABnew2eq1} when $n=3$ and $\omega\neq0$ ? In this paper we also
study this problem and obtain the following result.

\begin{theorem}\label{thmABABnew3}
Let $\omega, c, \lambda_{1}, \lambda_{2}, p_{1}, p_{2}$ be non-zero constants,   $q,Q$ be polynomials such that  $Q$ is not a constant and $q\not\equiv0$. If $f$ is a transcendental entire solution with finite order of
\begin{eqnarray}\label{thmABABnew2eq1addz2}
   f^{3}(z)+\omega f^{2}(z)f'(z)+q(z)e^{Q(z)}f(z+c)=p_{1}e^{\lambda_{1} z}+p_{2}e^{\lambda_{2} z},
\end{eqnarray}
 satisfying $N_{1)}(r,1/f)<(\kappa+o(1) T(r,f)$, where $0\leq \kappa<1$ and $N_{1)}(r,1/f)$ denotes the counting functions corresponding to simple  zeros of $f$,  then $\sigma(f)= \deg Q=1$.
\end{theorem}

The two examples below exhibit the sharpness of Theorem~\ref{thmABABnew3}.

\begin{example}\label{example3}
$f_{2}(z)=e^{z}$ is an transcendental entire solution of the nonlinear differential-difference equation
\begin{eqnarray*}
  f^{3}+f^{2}f'+\frac{1}{2}e^{-4z}f(z+\log 2)=2e^{3z}+e^{-3z}.
\end{eqnarray*}
Here $\omega=1\neq0$, $Q=-4z$,  $N_{1)}(r,1/f_{2})=0$ from the fact that $0$ is a Picard exceptional value of
$f_{2}$, and $\sigma(f_{2})= 1=\deg Q$.
\end{example}

\begin{example}\label{example4}
$f_{3}(z)=e^{2z}-e^{z}$ is an transcendental entire solution of the nonlinear differential-difference equation
\begin{eqnarray*}
  f^{3}-f^{2}f'-\frac{1}{5}e^{3z}f(z+\log 5)=-e^{6z}-3e^{5z}.
\end{eqnarray*}
Here $\omega=-1\neq0$, $Q=3z$, $N_{1)}(r,1/f_{3})=N(r,1/f_{3})=r/\pi+o(r)<
2r/\pi+o(r)=T(r,f_{3})$ by using the following Lemma~\ref{lemma24},
and $\sigma(f_{3})= 1=\deg Q$.
\end{example}

\section{Preliminary Lemmas}

The following two lemmas play important roles in
uniqueness problems of meromorphic functions.

\begin{lemma}[\cite{Yi1}]\label{lemma21}
 Let $f_{j}(z)\, (j=1,\ldots,n)\, (n\geq2)$ be meromorphic
functions, and let $g_{j}(z)\, (j=1,\ldots,n)$ be entire functions satisfying
\begin{itemize}
\item[(i)] $\sum _{j=1}^{n}f_{j}(z)e^{g_{j}(z)}\equiv 0;$
\item[(ii)]when $1\leq j < k \leq n,$ then $g_{j}(z)-g_{k}(z)$ is not a constant;
\item[(iii)] when $1\leq j \leq n, 1\leq h < k \leq n$, then
$$T(r,f_{j})=o \{T(r,e^{g_{h}-g_{k}})\} \quad (r\to \infty, r\not\in E),$$
where $E\subset(1,\infty)$ is of finite linear measure or logarithmic measure.
\end{itemize}
Then, $f_{j}(z)\equiv 0\, (j=1,\ldots ,n)$.
\end{lemma}

\begin{lemma}[\cite{Yi1}]\label{lemma27}
Let $f_{j}(z),\, j=1,2,3$ be  meromorphic functions and
$f_{1}(z)$ is not a constant. If
\begin{equation*}
    \sum_{j=1}^{3}f_{j}(z)\equiv 1,
  \end{equation*}
and
\begin{equation*}
    \sum_{j=1}^{3}N\left(r,\frac{1}{f_{j}}\right)
    +2  \sum_{j=1}^{3}\overline{N}(r,f_{j})
    <(\lambda+o(1))T(r),\quad r\in I,
  \end{equation*}
where $\lambda<1$, $T(r)=\max_{1\leq j\leq3}\{T(r,f_{j})\}$ and
$I$ represents a set of $r\in (0,\infty)$ with infinite linear measure. Then $f_{2}\equiv 1$ or $f_{3}\equiv 1$.
\end{lemma}

The difference analogues of Logarithmic Derivative Lemma (see \cite{chiang-Feng2009, Halburd-Korhonen20061, Halburd-Korhonen20062,Halburd-Korhonen2014, Korhonen}) play important roles in the
study of complex difference equations. The following version is a special case of \cite[Lemma 2.2]{Korhonen}.

\begin{lemma}\cite{Korhonen})\label{lemma22}
Let $f$ be a non-constant meromorphic function, let $c,h$ be two complex numbers such that $c\neq h$. If he hyper-order of $T(r,f)$ i.e.
$\sigma_{2}(f)<1$, then
\begin{equation*}
   m\left(r,\frac{f(z+h)}{f(z+c)} \right)=S(r,f)
 \end{equation*}
\end{lemma}
\noindent for all $r$ outside of a set of finite logarithmic measure.

The following lemma, which is a special case of \cite[Theorem 3.1]{Korhonen}, gives a relationship for the Nevanlinna
characteristics of meromorphic function with its shift.

\begin{lemma} \cite{Korhonen} \label{lemma23}
Let $f(z)$ be a  meromorphic function with the hyper-order less that
one, and $c\in \mathbb{C}\setminus \{0\}$. Then we have
 \begin{equation*}
   T(r,f(z+c))=T(r,f(z))+ S(r, f).
 \end{equation*}
\end{lemma}

The following lemma gives the Nevanlinna characteristic and counting functions of a exponential
polynomial. For convenience of the readers,  let's first recall some definitions and notations concerning  exponential polynomial of the form
\begin{eqnarray}\label{defeq1}
  f(z)=P_{1}(z)e^{Q_{1}(z)}+\cdots+P_{k}(z)e^{Q_{k}(z)},
\end{eqnarray}
where $P_{j}$ and  $Q_{j}$ are polynomials in $z$ for $1\leq j\leq k$.
Following Steinmetz \cite{Steinmetz}, \eqref{defeq1}
can be written in the normalized form
\begin{eqnarray}\label{defeq2}
  f(z)=H_{0}(z)+H_{1}(z)e^{\omega_{1}z^{q}}+\cdots+
H_{m}(z)e^{\omega_{m}z^{q}},
\end{eqnarray}
 where $H_{j}$ are either exponential polynomials of order $<q$ or ordinary polynomials in $z$, the leading coefficients $\omega_{j}$ are pairwise distinct, and $m\leq k$. In addition, the convex hull of a finite set $W\subset \mathbb{C}$, denoted by $co(W)$, is the intersection of finitely many closed half-planes each including $W$, and hence  $co(W)$ is either a compact polygon or a line segment. We denote the circumference of $co(W)$ by $C(co(W))$.
Concerning the exponential polynomial $f(z)$ in  \eqref{defeq2}, we denote
$W=\{\overline{\omega_{1}},\ldots,\overline{\omega_{m}}\}$ and $W_{0}=W\cup \{0\}$.

\begin{lemma}\cite{Steinmetz}\label{lemma24}
Let $f(z)$ be given by \eqref{defeq2}. Then
\begin{equation*}
   T(r,f)=C(co(W_{0}))\frac{r^{q}}{2\pi}+o(r^{q}).
 \end{equation*}
 If $H_{0}(z)\not\equiv 0$, then
 \begin{eqnarray*}
   m\left(r, \frac{1}{f}\right)=o(r^{q}),
 \end{eqnarray*}
 while if $H_{0}(z)\equiv0$, then
 \begin{eqnarray*}
    N\left(r, \frac{1}{f}\right)=C(co(W))\frac{r^{q}}{2\pi} + o(r^{q}).
 \end{eqnarray*}
\end{lemma}

The following lemma is a revised version of \cite[Lemma 2.4.2]{Laine}.
\begin{lemma}\label{lemma25}
Let $f(z)$ be a transcendental meromorphic solution of the equation:
 \begin{equation*}
   f^{n}P(z,f)=Q(z,f),
 \end{equation*}
where $P(z,f)$ and $Q(z,f)$ are polynomials in $f$ and
its derivatives with meromorhphic coefficients, say $\{a_{\lambda}|\lambda\in I\}$,
$n$ be a positive integer. If the total
degree of $Q(z,f)$ as a polynomial in $f$ and its derivatives
is at most $n$, then
\begin{eqnarray*}
  m(r,P(z,f))\leq \sum_{\lambda\in I}m(r,a_{\lambda})+S(r,f).
\end{eqnarray*}

\end{lemma}

\begin{lemma}(the Hadamard factorization theorem \cite[Theorem 2.7]{Yi1} or \cite[Theorem 1.9]{Conway}) \label{lemma26}
Let $f$ be a meromorphic function of finite order $\sigma(f)$. Write
\begin{eqnarray*}
  f(z)=c_{k}z^{k}+c_{k+1}z^{k+1}+\cdots \; (c_{k}\neq 0)
\end{eqnarray*}
near $z=0$ and let $\{a_{1}, a_{2}, \ldots\}$ and $\{b_{1}, b_{2}, \ldots\}$ be the zeros and poles of $f$ in $\mathbb{C}\backslash\{0\}$, respectively. Then
\begin{eqnarray*}
  f(z)=z^{k}e^{Q(z)}\frac{P_{1}(z)}{P_{2}(z)},
\end{eqnarray*}
where $P_{1}(z)$ and $P_{2}(z)$ are the canonical products of $f$ formed with
the non-null zeros and poles of $f(z)$, respectively, and $Q(z)$ is a polynomial
of degree $\leq \sigma(f)$.
\end{lemma}
\begin{remark} \label{remard21}
A well known fact about Lemma~\ref{lemma26} asserts that
$\lambda(f)=\lambda(z^{k}P_{1})=\sigma(z^{k}P_{1})\leq \sigma(f)$,  $\lambda(1/f)=\lambda(P_{2})=\sigma(P_{2})\leq \sigma(f)$ if $k\geq 0$; and
$\lambda(f)=\lambda(P_{1})=\sigma(P_{1})\leq \sigma(f)$,
 $\lambda(1/f)=\lambda(z^{-k}P_{2})=\sigma(z^{-k}P_{2})\leq \sigma(f)$ if $k<0$. So we have
$\sigma(f)=\sigma(e^{Q})$ when $\max\{\lambda(f), \lambda(1/f)\}<\sigma(f)$.
\end{remark}

By combining \cite[Theorem1.42]{Yi1} with
\cite[Theorem1.44]{Yi1}, we have the following lemma.

\begin{lemma}[\cite{Yi1}]\label{lemma28}
Let $f(z)$ be a non-constant meromorphic function in the complex plane. If
$0,\infty$ are Picard exceptional values of $f(z)$, then $f(z)=e^{h(z)}$, where $h(z)$ is a non-constant entire function.
Moreover,  $f(z)$ is of normal growth, and
\begin{itemize}
\item[(i)] if $h(z)$ is a polynomial of degree $p$, then $\sigma(f)=p$;
\item[(ii)] if $h(z)$ is a transcendental entire function, then  $\sigma(f)=\infty$.
\end{itemize}
\end{lemma}

The following lemma gives a relationship between the growth order of a function and its derivative.
\begin{lemma}[\cite{Yi1}]\label{lemma30}
Suppose that $f(z)$ is meromorphic in the complex plane and
$n$ is a positive integer. Then
$f(z)$ and $f^{(n)}(z)$ have the same order.
\end{lemma}

\section{Proof of Theorem~\ref{thmABABnew1}.}
Let $f$ be a  transcendental entire solution with finite order of equation \eqref{aaa111} satisfying $\lambda(f)<\sigma(f)$. Then, by Lemma~\ref{lemma26} and Remark~\ref{remard21},  we can factorize $f(z)$ as
\begin{eqnarray}\label{thmA1pfeqnnm1}
 f(z)=d(z)e^{h(z)},
\end{eqnarray}
where $h$ is a polynomial with $\deg h=\sigma(f)$, $d$
 is the canonical products formed by zeros of $f$ with $\sigma(d)=\lambda(f)<\sigma(f)$.
Obviously, $h$ is a non-constant polynomial, otherwise we'll have $\sigma(f)=\sigma(d)=\lambda(f)$,  a contradiction. Thus we have that
$\deg h\geq 1$. Let $\deg h=m\, (\geq 1)$, and $h(z)=a_{m}z^{m}+a_{m-1}z^{m-1}+\cdots$, where $a_{m}\neq 0$.

We rewrite \eqref{aaa111}  as
\begin{eqnarray}\label{aaa113}
 f^{2}+\widetilde{\omega} f f'+qe^{Q}f_{c}=ue^{v},
\end{eqnarray}
where $f_{c}=f(z+c)$, for simplicity.

Obviously,
we have $\sigma(f_{c})=\sigma(f)=\sigma(f')$ by Lemma~\ref{lemma23}
and Lemma~\ref{lemma30}. So from \eqref{aaa113}, by the order property,
  we get
\begin{eqnarray}\label{addaaa113}
\deg v=\sigma(ue^{v})&\leq& \max \{\sigma(f')=\sigma(f)=\sigma(f_{c}), \sigma(e^{Q}), \sigma(q) \}\nonumber \\
&=&\max \{\deg h, \deg Q \} .
\end{eqnarray}

By substituting \eqref{thmA1pfeqnnm1} into
\eqref{aaa113},  we get
\begin{eqnarray}\label{add111}
  \left(d^{2}+\widetilde{\omega} d\cdot (d'+dh')\right)e^{2h}+qd_{c}e^{Q+h_{c}}=ue^{v}.
\end{eqnarray}

{\bf Case 1.} $\sigma(f)> \deg Q$,  then we have $\deg h> \deg Q\geq 1$, and
$\deg v\leq \deg h$ from \eqref{addaaa113}.

{\bf Subcase 1.1.} $\deg h>\deg v$.
From \eqref{add111} we have
\begin{eqnarray}\label{addaddaddzzz111}
 \left(d^{2}+\widetilde{\omega} d\cdot (d'+dh')\right)e^{h_{1}}
 e^{2a_{m}z^{m}}
 +qd_{c}e^{h_{2}}e^{a_{m}z^{m}}
  =ue^{v},
\end{eqnarray}
where $h_{1}=2a_{m-1}z^{m-1}+\cdots$ and  $h_{2}=Q+(a_{m}mc+a_{m-1})z^{m-1}+\cdots$ are all polynomials with degree at most $m-1$. So, combining with $\sigma(d')=\sigma(d)=\sigma(d_{c})<m$,
 by using Lemma~\ref{lemma21} to \eqref{addaddaddzzz111}, we have
\begin{eqnarray*}
qd_{c} \equiv 0,
\end{eqnarray*}
 which
yields a contradiction. Thus $\deg h>\deg v$ can not holds.

{\bf Subcase 1.2.} $\deg h=\deg v$. Let $v(z)=v_{m}z^{m}+v_{m-1}z^{m-1}+\cdots $, where $v_{m}\neq 0$. From \eqref{add111} we have
\begin{eqnarray}\label{theoremeq222}
 \left(d^{2}+\widetilde{\omega} d\cdot (d'+dh')\right)e^{h_{1}}
 e^{2a_{m}z^{m}}
 +qd_{c}e^{h_{2}}e^{a_{m}z^{m}}
  =ue^{h_{3}}e^{v_{m}z^{m}}.
\end{eqnarray}
where $h_{3}=v_{m-1}z^{m-1}+\cdots$ is a polynomial with degree at most $m-1$, $h_{1}$ and $h_{2}$ are defined as in Subcase 1.1.

If  $v_{m}\neq 2a_{m}$ and $v_{m}\neq a_{m}$, combining with $\sigma(d')=\sigma(d)=\sigma(d_{c})<m$,  by using Lemma~\ref{lemma21} to \eqref{theoremeq222},
we get $u\equiv 0$, a contradiction.

If $v_{m}=2a_{m}$, then \eqref{theoremeq222} can be reduced to
\begin{eqnarray*}
 \left( \left(d^{2}+\widetilde{\omega} d\cdot (d'+dh')\right)e^{h_{1}}-ue^{h_{3}}
 \right)e^{2a_{m}z^{m}}
 +qd_{c}e^{h_{2}}e^{a_{m}z^{m}}
  =0.
\end{eqnarray*}
 Thus, by using Lemma~\ref{lemma21}, we have $qd_{c}\equiv 0$, a contradiction.

If $v_{m}=a_{m}$, then \eqref{theoremeq222} can be reduced to
\begin{eqnarray*}
  \left(d^{2}+\widetilde{\omega} d\cdot (d'+dh')\right)e^{h_{1}}
 e^{2a_{m}z^{m}}
 +(qd_{c}e^{h_{2}}-ue^{h_{3}})e^{a_{m}z^{m}}
  =0.
\end{eqnarray*}
Similarly as above, by Lemma~\ref{lemma21}, we get
\begin{eqnarray*}
 d^{2}+\widetilde{\omega} d\cdot (d'+dh')\equiv 0.
\end{eqnarray*}
This gives that
\begin{eqnarray*}
 d=c_{1}e^{-\frac{1}{\widetilde{\omega}}z-h}, \quad  c_{1}\in \mathbb{C}\setminus \{0\}.
\end{eqnarray*}
Since  $\deg h> \deg Q\geq 1$, so we have $\sigma(d)=\deg h=\sigma(f)$, which
contradicts with the assumption that $\sigma(d)<\sigma(f)$.

{\bf Case 2.} $\sigma(f)=\deg Q$. Let $Q(z)=b_{m}z^{m}+b_{m-1}z^{m-1}+\cdots$,
where $b_{m}\neq0$.

Suppose $\deg v< \deg h=\deg Q$. From \eqref{add111} we have
\begin{eqnarray}\label{theoremeq111}
  \left(d^{2}+\widetilde{\omega} d\cdot (d'+dh')\right)e^{h_{1}}
 e^{2a_{m}z^{m}}
 +(qd_{c}e^{\widetilde{h_{2}}})e^{(a_{m}+b_{m})z^{m}}=ue^{v},
\end{eqnarray}
where $\widetilde{h_{2}}=(a_{m}mc+a_{m-1}+b_{m-1})z^{m-1}+\cdots$ is a
polynomial with degree at most $m-1$, and $h_{1}$ is defined as in Subcase 1.1.

If $b_{m}\neq \pm a_{m}$, combining with $\sigma(d')=\sigma(d)=\sigma(d_{c})<m$,  by using Lemma~\ref{lemma21} to \eqref{theoremeq111}, we get $u\equiv0$, which yields a contradiction.

If $b_{m}=a_{m}$, then \eqref{theoremeq111}
can be reduced to
\begin{eqnarray*}
\left( \left(d^{2}+\widetilde{\omega} d\cdot (d'+dh')\right)e^{h_{1}}
 + qd_{c}e^{\widetilde{h_{2}}}  \right)e^{2a_{m}z^{m}}
 =ue^{v}.
\end{eqnarray*}
Thus, by using Lemma~\ref{lemma21}, we have $u\equiv 0$, a contradiction.

If $b_{m}=-a_{m}$, then \eqref{theoremeq111}
can be reduced to
\begin{eqnarray*}
  \left(d^{2}+\widetilde{\omega} d\cdot (d'+dh')\right)e^{h_{1}}
 e^{2a_{m}z^{m}}
 =ue^{v}-qd_{c}e^{\widetilde{h_{2}}}.
\end{eqnarray*}
So by  Lemma~\ref{lemma21}, we get
\begin{eqnarray*}
  d^{2}+\widetilde{\omega} d\cdot (d'+dh')\equiv 0.
\end{eqnarray*}
This gives that
\begin{eqnarray*}
 d=c_{2}e^{-\frac{1}{\widetilde{\omega}}z-h},\; c_{2}\in\mathbb{C}\setminus\{0\}.
\end{eqnarray*}
Thus by $\deg h>\deg v\geq 1$, we have $\sigma(d)=\deg h=\sigma(f)$,
which contradicts with the fact that $\sigma(d)<\sigma(f)$.

Therefore, we have $\deg v= \deg h=\deg Q$ from \eqref{addaaa113}, which implies that Conclusion (2) holds.

{\bf Case 3.} $\sigma(f)< \deg Q$, then we have $T(r,f)=S(r,e^{Q})$.
Thus we get $T(r,f')=S(r,e^{Q})$ from Milloux's theorem and $T(r,f_{c})=S(r,e^{Q})$ from lemma~\ref{lemma23}. Therefore,
by \eqref{aaa113}, we have
\begin{eqnarray*}
 T(r,e^{Q})+S(r,e^{Q})&=&T(r,f^{2}+\widetilde{\omega} ff'+qf_{c}e^{Q})\\
 &=&T(r,ue^{v})=T(r,e^{v})+S(r,e^{v}).
\end{eqnarray*}

Therefore,
\begin{eqnarray*}\label{aaa112}
 \deg Q = \deg v.
\end{eqnarray*}
 Differentiating \eqref{aaa113} yields
\begin{eqnarray}\label{aaa114}
  2ff'+\widetilde{\omega} (f')^{2}+\widetilde{\omega} f f''+Ae^{Q}=(u'+uv')e^{v},
\end{eqnarray}
with $A=q'f_{c}+qf_{c}'+qf_{c}Q'$.

Eliminating $e^{v}$ from \eqref{aaa113} and \eqref{aaa114} to get
\begin{eqnarray}\label{aaa115}
 B_{1}e^{Q}+B_{2}=0,
\end{eqnarray}
where
$$B_{1}=uA-qf_{c}(u'+uv'),$$
$$B_{2}=u[2ff'+\widetilde{\omega} (f')^{2}+\widetilde{\omega} ff'']-(f^{2}+\widetilde{\omega} ff')(u'+uv').$$

Noticing that $\sigma(f_{c})=\sigma(f)< \deg Q$, and $\sigma(f'')=\sigma(f')=\sigma(f)< \deg Q$ from
Lemma~\ref{lemma30}, thus
by Lemma~\ref{lemma21}, we get $B_{1}\equiv B_{2}\equiv 0$. It follows from $B_{1}\equiv 0$ that
\begin{eqnarray*}
 \frac{q'}{q}+\frac{f_{c}'}{f_{c}}+Q' =\frac{u'}{u}+v',
\end{eqnarray*}
by integrating, we have $qf_{c}e^{Q}=c_{3}ue^{v}$, where $c_{3}$ is
a non-zero constant.

If $c_{3}=1$, by substituting $qf_{c}e^{Q}=ue^{v}$ into \eqref{aaa113}, we see that $f^{2}+\widetilde{\omega} f f'=0$. Thus we can easily get $f=c_{4}e^{-z/\widetilde{\omega}}$,
where $c_{4}\in \mathbb{C}\setminus \{0\}$, which implies that Conclusion (1) holds.

If $c_{3}\neq1$, we have $f=c_{3}u_{-c}/q_{-c}e^{v_{-c}-Q_{-c}}$. By substituting it into \eqref{aaa113}, we get
\begin{eqnarray*}
 \frac{c_{3}u_{-c}}{q_{-c}}  \left(\frac{c_{3}u_{-c}}{q_{-c}}
  +\widetilde{\omega}  \left( \left( \frac{c_{3}u_{-c}}{q_{-c}}\right)'  +\frac{c_{3}u_{-c}}{q_{-c}} (v_{-c}-Q_{-c}) ' \right) \right)e^{2(v_{-c}-Q_{-c})}
  =(1-c_{3})ue^{v}
\end{eqnarray*}
Since  $1\leq\deg h=\sigma(f)<\deg Q=\deg v$ and $\sigma(f)=\deg(v_{-c}-Q_{-c})$, so by
 Lemma~\ref{lemma21} we can easily deduce a contradiction
 by the fact that   $c_{3}\neq1$ and $u\not\equiv 0$.

\section{Proof of Theorem~\ref{thmABABnew2}.}

Suppose that $f$ is a transcendental entire solution with finite order of equation \eqref{thmABABnew2eq1}. We rewrite \eqref{thmABABnew2eq1}  as
\begin{eqnarray}\label{zzaaa113}
 f^{n}+\omega f^{n-1}f'+qe^{Q}f_{c}=p_{1}e^{\lambda_{1} z}+p_{2}e^{\lambda_{2} z},
\end{eqnarray}
where $f_{c}=f(z+c)$, for short. From Lemma~\ref{lemma23}, we have $\sigma(f)=\sigma(f_{c})$.

By differentiating both sides of \eqref{zzaaa113}, we have
\begin{eqnarray}\label{zzaaa113z1}
nf^{n-1}f'+\omega(n-1)f^{n-2}(f')^{2}+
\omega f^{n-1}f''+A_{1}e^{Q}
=p_{1}\lambda_{1}e^{\lambda_{1} z}+p_{2}\lambda_{2}e^{\lambda_{2} z},
\end{eqnarray}
where $A_{1}=q'f_{c}+qf_{c}'+qf_{c}Q'$.

By eliminating $e^{\lambda_{2} z}$ from equations \eqref{zzaaa113}
and \eqref{zzaaa113z1}, we get
\begin{eqnarray}\label{zzaaa113z2}
&& \lambda_{2}f^{n}+(\lambda_{2}\omega-n) f^{n-1}f'
-\omega(n-1)f^{n-2}(f')^{2}-
\omega f^{n-1}f''+A_{2}e^{Q}\nonumber\\
&&=p_{1}(\lambda_{2}-\lambda_{1})e^{\lambda_{1} z},
\end{eqnarray}
where $A_{2}=\lambda_{2}qf_{c}-A_{1}$.

By differentiating \eqref{zzaaa113z2}, we have
\begin{eqnarray}\label{zzaaa113z3}
&&\lambda_{2}nf^{n-1}f'+(\lambda_{2}\omega-n)\left[(n-1)f^{n-2}(f')^{2}+ f^{n-1}f''\right]\nonumber\\
&&-\omega(n-1)\left[(n-2)f^{n-3}(f')^{3}+ f^{n-2}2f'f'' \right]
 -\omega (n-1)f^{n-2}f'f''-\omega f^{n-1}f'''
\nonumber\\
&&+(A_{2}'+A_{2}Q')e^{Q}
=p_{1}(\lambda_{2}-\lambda_{1})\lambda_{1}e^{\lambda_{1} z}.
\end{eqnarray}

By eliminating $e^{\lambda_{1} z}$ from equations \eqref{zzaaa113z2}
and \eqref{zzaaa113z3}, we obtain
\begin{eqnarray}\label{zzaaa113z4}
&&\lambda_{1}\lambda_{2}f^{n}+(\lambda_{1}\lambda_{2}\omega-n\lambda_{1}-\lambda_{2}n) f^{n-1}f'
-(n-1)\left[\lambda_{1}\omega +\lambda_{2}\omega-n \right]f^{n-2}(f')^{2}\nonumber\\
&&-(\lambda_{1}\omega+\lambda_{2}\omega-n) f^{n-1}f''+\omega(n-1)(n-2)f^{n-3}(f')^{3}
 +3\omega (n-1)f^{n-2}f'f''\nonumber\\
&&+\omega f^{n-1}f'''+(\lambda_{1}A_{2}-A_{2}'-A_{2}Q')e^{Q}
=0.
\end{eqnarray}

{\bf Case 1.}  $\sigma(f)<1$.  By using  logarithmic derivative
lemma, Lemma~\ref{lemma22} and Lemma~\ref{lemma24}, from \eqref{zzaaa113}  we get
\begin{eqnarray*}
  T\left(r,e^{Q}\right)&=&m\left(r,e^{Q}\right)=m\left(r, \frac{p_{1}e^{\lambda_{1}z}+p_{2}e^{\lambda_{2}z}-f^{n}-\omega f^{n-1}f'}{qf_{c}} \right)\\
  &\leq& m\left(r,\frac{f}{qf_{c}}\right)+m\left(r,\frac{1}{f} \right)+m\left(r,p_{1}e^{\lambda_{1}z}+p_{2}e^{\lambda_{2}z} \right)\\
  &&+m\left(r,\frac{f^{n}+\omega f^{n-1}f'}{f^{n}} \right)+m(r,f^{n})
  +O(1)\\
  &\leq& (n+1)T(r,f)+C(co(W_{0}))\frac{r}{2\pi}+o(r)+S(r,f)\\
  &\leq& C(co(W_{0}))\frac{r}{2\pi}+o(r),
\end{eqnarray*}
where $W_{0}=\{0,\overline{\lambda_{1}}, \overline{\lambda_{2}}\}$, thus we
have $\deg Q \leq 1$, and noting that $\deg Q\geq 1$, we know that $\deg Q=1$. We set $Q=az+b$, $a\in \mathbb{C}\setminus \{0\}$, $b\in \mathbb{C}$.

Thus, by using Lemma~\ref{lemma21} to \eqref{zzaaa113z4}, we have
\begin{eqnarray}\label{theorem2eq0}
\lambda_{1}A_{2}-A_{2}'-A_{2}Q'=(\lambda_{1}-a)A_{2}-A_{2}'\equiv 0.
\end{eqnarray}

{\bf Subcase 1.1.} $A_{2}\equiv 0$. That is $\lambda_{2}qf_{c}-q'f_{c}-qf_{c}'-qf_{c}a \equiv 0.$
This gives that
\begin{eqnarray*}
 \lambda_{2}-\frac{q'}{q}-\frac{f_{c}'}{f_{c}}-a \equiv 0.
\end{eqnarray*}
By integrating, we have
\begin{eqnarray*}
qf_{c}=c_{1}e^{(\lambda_{2}-a)z}, \quad c_{1}\in \mathbb{C}\setminus \{0\}.
\end{eqnarray*}
So we have $a=\lambda_{2}$, and $f_{c}=c_{1}/q$.
Otherwise, if $a\neq\lambda_{2}$, then we get
$\sigma(f)=\sigma(f_{c})=1$, which contradicts with
our assumption that $\sigma(f)<1$. So we have that
$f(z)=c_{1}/q(z-c)$ is a rational function, which contradicts with the
assumption that $f$ is transcendental.

{\bf Subcase 1.2.} $A_{2}\not\equiv 0$. From \eqref{theorem2eq0}, we
get
\begin{eqnarray*}
A_{2}
=c_{2}e^{(\lambda_{1}-a)z}, \quad c_{2}\in \mathbb{C}\setminus \{0\}.
\end{eqnarray*}
This gives that
\begin{eqnarray}\label{theorem2eq1}
(qf_{c})'+(a-\lambda_{2})(qf_{c})
=-c_{2}e^{(\lambda_{1}-a)z}, \quad c_{2}\in \mathbb{C}\setminus \{0\}.
\end{eqnarray}
Since $\lambda_{1}\neq \lambda_{2}$, we discuss the following three subcases:

{\bf Subcase 1.2.1.} $a=\lambda_{2}$. Then \eqref{theorem2eq1} reduces to $(qf_{c})'=-c_{2}e^{(\lambda_{1}-\lambda_{2})z}$. Thus,
we have $qf_{c}=\frac{c_{2}}{\lambda_{2}-\lambda_{1}}
e^{(\lambda_{1}-\lambda_{2})z}+c_{3}$, where $c_{3}\in \mathbb{C}$. Therefore, by  Lemma~\ref{lemma24} we have
\begin{eqnarray*}
T(r,qf_{c})&=&
T\left(r,\frac{c_{2}}{\lambda_{2}-\lambda_{1}}
e^{(\lambda_{1}-\lambda_{2})z} +c_{3} \right)\\
&=& C(co(W_{1}))\frac{r}{2\pi}+o(r), \quad W_{1}=\{0,\overline{\lambda_{1}-\lambda_{2}}\}.
\end{eqnarray*}
Combining with Lemma~\ref{lemma23} and the fact that $f$ is transcendental, we get
\begin{eqnarray*}
C(co(W_{1}))\frac{r}{2\pi}+o(r)=T(r,qf_{c})&\leq& T(r,f_{c})+T(r,q)+O(1)\\
&=& T(r,f)+S(r,f),
\end{eqnarray*}
which contradicts with the assumption that $\sigma(f)<1$.

{\bf Subcase 1.2.2.} $a=\lambda_{1}$. Then \eqref{theorem2eq1} reduces to $(qf_{c})'+(\lambda_{1}-\lambda_{2})(qf_{c})
=-c_{2}$. Thus, we have $qf_{c}=\frac{c_{2}}{\lambda_{2}-\lambda_{1}}
+c_{4}e^{(\lambda_{2}-\lambda_{1})z}$, where $c_{4}\in \mathbb{C}$.
We assert that $c_{4}\neq 0$. Otherwise, if $c_{4}=0$, we get
that $f(z)=\frac{c_{2}}{\lambda_{2}-\lambda_{1}}\frac{1}{q(z-c)}$ is rational, which contradicts with the assumption that $f$ is transcendental. Therefore, since $c_{4}\neq 0$ and $\lambda_{1}\neq \lambda_{2}$, similarly as in Subcase 1.2.1, by combining with Lemma~\ref{lemma23}, Lemma~\ref{lemma24}, and the assumption that $\sigma(f)<1$, we can also get a contradiction.

{\bf Subcase 1.2.3.} $a\neq\lambda_{1}$ and $a\neq\lambda_{2}$.
Then by \eqref{theorem2eq1}, we get that
\begin{eqnarray*}
qf_{c}=\frac{c_{2}}{\lambda_{2}-\lambda_{1}}e^{(\lambda_{1}-a)z}
+c_{5}e^{(\lambda_{2}-a)z}, \quad   c_{5}\in \mathbb{C}.
\end{eqnarray*}
Therefore, since $c_{2}\neq 0$, $a\neq\lambda_{1}$,   and
$\lambda_{1}\neq \lambda_{2}$, similarly as in Subcase 1.2.1, by combining with Lemma~\ref{lemma23}, Lemma~\ref{lemma24}, and the assumption that $\sigma(f)<1$, we also get a contradiction.

{\bf Case 2.} $\sigma(f)>1$. Denote $P=p_{1}e^{\lambda_{1} z}+p_{2}e^{\lambda_{2} z}$,  then we have $\sigma(P)=1$ by Lemma~\ref{lemma24},
and equation \eqref{zzaaa113} can be rewritten as:
\begin{eqnarray}\label{theoremeqz01}
 f^{n}+\omega f^{n-1}f'+(qf_{c})e^{Q}=P.
\end{eqnarray}
Differentiating \eqref{theoremeqz01} yields
\begin{eqnarray}\label{theoremeqz02}
 nf^{n-1}f'+\omega (n-1)f^{n-2}(f')^{2}+ \omega f^{n-1}f''
 +Le^{Q}=P',
\end{eqnarray}
where $L=(qf_{c})'+Q'(qf_{c})$.

{\bf Subcase 2.1.} $\omega\neq 0$ and $n\geq4$.
Eliminating $e^{Q}$ from \eqref{theoremeqz01} and \eqref{theoremeqz02}, we have
\begin{eqnarray}\label{theoremeqz03}
 f^{n-2}H=PL-P'(qf_{c}),
\end{eqnarray}
where
\begin{eqnarray*}
H=Lf^{2}+(\omega L-nqf_{c}) ff'-(n-1)\omega qf_{c}(f')^{2}-\omega qf_{c} f f''.
\end{eqnarray*}

{\bf Subcase 2.1.1.} $H\not\equiv 0$. From the assumptions that $f$ is entire and $q,Q$ are polynomials, we have that both $H$ and $fH$ are entire functions.  We rewrite
$PL-P'(qf_{c})$ as
\begin{eqnarray*}
 && \quad P\left[(qf_{c})'+Q'(qf_{c})\right]-P'(qf_{c})\\
  &&  =Pq\frac{(qf_{c})'}{qf_{c}}\frac{f_{c}}{f}\cdot f+(PQ'-P')q\frac{f_{c}}{f}\cdot f\\
  && =Pq\left(\frac{q'}{q}+\frac{f_{c}'}{f_{c}}  \right)\frac{f_{c}}{f}\cdot f+(PQ'-P')q\frac{f_{c}}{f}\cdot f,
\end{eqnarray*}
and rewrite $H$ as
\begin{eqnarray*}
 &&q \left(\frac{q'}{q}+\frac{f_{c}'}{f_{c}}  +Q'\right)\frac{f_{c}}{f}\cdot f^{3}
 +q\left(\omega  \left(\frac{q'}{q}+\frac{f_{c}'}{f_{c}} + Q' \right) -n\right)\frac{f_{c}}{f}\cdot f^{2}f'\\
 &&\, -(n-1)\omega q\frac{f_{c}}{f}\cdot f(f')^{2}-\omega q\frac{f_{c}}{f}\cdot f^{2} f'',
\end{eqnarray*}
thus both $PL-P'(qf_{c})$ and $H$ are all differential polynomials with meromorphic coefficients. By logarithmic derivative lemma and Lemma~\ref{lemma23}, we have $m(r,f_{c}'/f_{c})=S(r,f_{c})=S(r,f)$; by Lemma~\ref{lemma22},
we have $m(r,f_{c}/f)=S(r,f)$; and by Lemma~\ref{lemma24}, we have
$m(r,P)=O(r)$. Therefore, by using
Lemma~\ref{lemma25} to \eqref{theoremeqz03} and the fact that $n\geq 4$ when
$\omega\neq 0$, we obtain that
\begin{eqnarray*}
T(r,H)=m(r,H)=S(r,f)+O(r).
\end{eqnarray*}
and
\begin{eqnarray*}
T(r,fH)=m(r,fH)=S(r,f)+O(r).
\end{eqnarray*}
Thus, by $H\not\equiv 0$ we have
\begin{eqnarray*}
T(r,f)\leq T(r,fH)+T\left(r, \frac{1}{H}\right)=S(r,f)+O(r),
\end{eqnarray*}
which contradicts with the assumption that $\sigma(f)>1$.

{\bf Subcase 2.1.2.} $H\equiv 0$. Then from \eqref{theoremeqz03}, we have
\begin{eqnarray*}
  PL-P'(qf_{c})=P\left[(qf_{c})'+Q'(qf_{c})\right]-P'(qf_{c})\equiv 0.
\end{eqnarray*}
This gives that
\begin{eqnarray*}
  \frac{(qf_{c})'}{qf_{c}}+Q'-\frac{P'}{P}\equiv 0.
\end{eqnarray*}
By integration, we see that there exists a $c_{6}\in \mathbb{C}\setminus\{0\}$ such that
\begin{eqnarray}\label{theoremeqz05}
 qf_{c}=c_{6}Pe^{-Q}.
\end{eqnarray}
By substituting \eqref{theoremeqz05} into \eqref{theoremeqz01}, we get that
\begin{eqnarray}\label{theoremeqz06}
 f^{n}+\omega f^{n-1}f'=(1-c_{6})P.
\end{eqnarray}
If $c_{6}=1$, then from \eqref{theoremeqz06} we have $f+\omega f'=0$. By integration, we get that $f=c_{7}e^{-\frac{1}{\omega}z}$, $c_{7}\in \mathbb{C}\setminus\{0\}$, which contradicts with the assumption that $\sigma(f)>1$. So we have $c_{6}\neq1$. From \eqref{theoremeqz05}, we have
\begin{eqnarray}\label{theoremeqz07}
  f=c_{6}\frac{P_{-c}}{q_{-c}}e^{-Q_{-c}},
\end{eqnarray}
and $\deg Q=\deg Q_{-c}=\sigma(f)>1$ since $\sigma(P_{-c})=1$ by lemma~\ref{lemma24}.

By Substituting \eqref{theoremeqz07} into \eqref{theoremeqz06}, we have
\begin{eqnarray*}
\frac{ c_{6}^{n}}{1-c_{6}} \left(\left(\frac{P_{-c}}{q_{-c}}\right)^{n}
  +\omega \left(\frac{P_{-c}}{q_{-c}}\right)^{n-1}
  \left( \left(\frac{P_{-c}}{q_{-c}}\right)'
  +\frac{P_{-c}}{q_{-c}}(-Q_{-c})'  \right)  \right)e^{-nQ_{-c}}
  =P.
\end{eqnarray*}
Then from $\deg Q>1$, $\sigma(P_{-c})=\sigma(P)=1$ and Lemma~\ref{lemma21},
we get that  $P(z)\equiv0$, a contradiction.

{\bf Subcase 2.2.} $\omega=0$ and $n\geq 3$.
Eliminating $e^{Q}$ from \eqref{theoremeqz01} and \eqref{theoremeqz02}, we have
\begin{eqnarray*}
 f^{n-1}\left( Lf-nqf_{c} f' \right)=PL-P'(qf_{c}).
\end{eqnarray*}

{\bf Subcase 2.2.1.} $Lf-nqf_{c} f'\not\equiv 0$. Since $n\geq 3$ and
$\omega= 0$, similarly as in subcase 2.1.1, we have
\begin{eqnarray*}
T(r, Lf-nqf_{c} f' )=m(r,Lf-nqf_{c} f')=S(r,f)+O(r).
\end{eqnarray*}
and
\begin{eqnarray*}
T(r,f\left( Lf-nqf_{c} f' \right))=m(r,f\left( Lf-nqf_{c} f' \right))=S(r,f)+O(r).
\end{eqnarray*}
Thus, by $Lf-nqf_{c} f'\not\equiv 0$ we have
\begin{eqnarray*}
T(r,f)\leq T(r,f\left( Lf-nqf_{c} f' \right))+T\left(r, \frac{1}{Lf-nqf_{c} f'}\right)=S(r,f)+O(r),
\end{eqnarray*}
which contradicts with the assumption that $\sigma(f)>1$.

{\bf Subcase 2.2.2.} $Lf-nqf_{c} f'\equiv 0$. Then
\begin{eqnarray*}
  \frac{(qf_{c})'}{qf_{c}}+Q'-n \frac{f'}{f}\equiv 0.
\end{eqnarray*}
By integration, we see that there exists a $c_{8}\in \mathbb{C}\setminus \{0\}$  such that
\begin{eqnarray}\label{theoremeqz09}
 qf_{c}e^{Q}=c_{8}f^{n}.
\end{eqnarray}
Substituting \eqref{theoremeqz09} into \eqref{theoremeqz01} gives
\begin{eqnarray*}
 (1+c_{8})f^{n}=P.
\end{eqnarray*}
If $c_{8}\neq -1$, we have $nT(r,f)+S(r,f)=T(r,(1+c_{8})f^{n})=T(r,P)=O(r)$, which contradicts with the assumption that $\sigma(f)>1$. If $c_{8}=-1$, then $P=p_{1}e^{\lambda_{1} z}+p_{2}e^{\lambda_{2} z}\equiv 0$, a contradiction.

{\bf Case 3.}  $\sigma(f)=1$. By \eqref{zzaaa113}, Lemma~\ref{lemma22}, and logarithmic derivative lemma, we obtain
\begin{eqnarray*}
  T\left(r,e^{Q}\right)&=& m\left(r,e^{Q}\right)=
  m\left(r,\frac{p_{1}e^{\lambda_{1} z}+p_{2}e^{\lambda_{2} z}-f^{n}-\omega f^{n-1}f'}{qf_{c}} \right)\\
&\leq&  m\left(r,\frac{1}{qf_{c}} \right)+m\left(r,p_{1}e^{\lambda_{1} z}+p_{2}e^{\lambda_{2} z}\right)+m\left(r, f^{n}+\omega f^{n-1}f'\right)+O(1)\\
&\leq&  m\left(r,\frac{f}{f_{c}} \right)+m\left(r,\frac{1}{f}\right)
+m\left(r, \frac{f^{n}+\omega f^{n-1}f'}{f^{n}}\right)+m(r,f^{n})\\
&&\quad +T
\left(r,p_{1}e^{\lambda_{1} z}+p_{2}e^{\lambda_{2} z}\right)+O(\log r)\\
&\leq& (n+1)T(r,f)+T
\left(r,p_{1}e^{\lambda_{1} z}+p_{2}e^{\lambda_{2} z}\right)+S(r,f).
\end{eqnarray*}
Note that $\deg Q \geq 1$, then by combining with Lemma~\ref{lemma24}, we get
\begin{eqnarray*}
  1\leq \deg Q= \sigma\left(e^{Q}\right)\leq
  \max\{\sigma(f), \sigma\left(p_{1}e^{\lambda_{1} z}+p_{2}e^{\lambda_{2}z}\right)\}=1,
\end{eqnarray*}
that is $\sigma(f)=\deg Q=1$. Thus,  Conclusion (1) is proved.

Next, we prove Conclusion (2). Suppose that $f$ is a transcendental entire solution with finite order of
 equation \eqref{thmABABnew2eq1}  with $\lambda(f)<\sigma(f)$.
Then, by Lemma~\ref{lemma26} and Remark~\ref{remard21},  we can factorize $f(z)$ as
\begin{eqnarray}\label{theoremeqz10}
 f(z)=d(z)e^{h(z)},
\end{eqnarray}
where $h$ is a polynomial with $\deg h=\sigma(f)$, $d$
 is the canonical products formed by zeros of $f$ with $\sigma(d)=\lambda(f)<\sigma(f)$. Similarly as in the proof of Theorem~\ref{thmABABnew1}, we have $\sigma(f)=\deg h\geq1$.

By substituting \eqref{theoremeqz10} into \eqref{thmABABnew2eq1}, we
get
\begin{eqnarray}\label{theoremeqz11}
 d^{n-1}(d+\omega(d'+dh'))e^{nh}+qd_{c}e^{Q+h_{c}}=p_{1}e^{\lambda_{1}z}
 +p_{2}e^{\lambda_{2}z}.
\end{eqnarray}

Dividing both sides of \eqref{theoremeqz11} by $p_{2}e^{\lambda_{2}z}$, we obtain \begin{eqnarray}\label{theoremeqz12}
 f_{1}+f_{2}+f_{3}
 =1,
\end{eqnarray}
where
\begin{eqnarray*}
 f_{1}= -\frac{p_{1}}{p_{2}}e^{(\lambda_{1}-\lambda_{2})z},
  f_{2}=\frac{d^{n-1}(d+\omega(d'+dh'))}{p_{2}}e^{nh-\lambda_{2}z},
  f_{3}= \frac{qd_{c}}{p_{2}}e^{Q+h_{c}-\lambda_{2}z}.
\end{eqnarray*}

Obviously, $f_{1}$ is
not a constant since $\lambda_{1}\neq\lambda_{2}$. Next we discuss two cases: $\sigma(f)>1$ and $\sigma(f)=1$, respectively. We set
$T(r)=\max\{T(r, f_{1}), T(r, f_{2}), T(r, f_{3})\}$.

If $\sigma(f)>1$, then $d^{n-1}(d+\omega(d'+dh'))/p_{2}e^{-\lambda_{2}z}$ is a small function of $e^{h}$ since $\max\{\sigma(d')= \sigma(d),1\}<\deg h$. So we have $T(r)\geq T(r, f_{2})=nT(r,e^{h})+S(r,e^{h})$ . Thus, by Milloux's theorem and Lemma~\ref{lemma23}, we get
\begin{eqnarray*}
 \frac{ N\left(r,\frac{1}{f_{2}} \right)}{T(r)}&=&\frac{N\left(r,\frac{1}{d^{n-1}(d+\omega(d'+dh'))} \right)}{T(r)}
  \leq \frac{T\left(r, d^{n-1}(d+\omega(d'+dh')) \right)}{T(r)}
  \\
  &=&\frac{O(T(r,d))+O(\log r)}{T(r,e^{h})}\cdot\frac{T(r,e^{h})}{T(r)}
  \to 0,
\end{eqnarray*}
and
\begin{eqnarray*}
 \frac{N\left(r,\frac{1}{f_{3}} \right)}{T(r)} &=&
 \frac{N\left(r,\frac{1}{qd_{c}} \right)}{T(r)}\leq
 \frac{T(r,qd_{c})}{T(r)}\leq
 \frac{T(r,d)+S(r,d)+O(\log r)}{T(r,e^{h})}\cdot\frac{T(r,e^{h})}{T(r)}
  \to 0
\end{eqnarray*}
as $ r\to \infty$.

Therefore, by using
Lemma~\ref{lemma27}, we can deduce that $f_{2}\equiv 1$ or $f_{3}\equiv 1$.

If $f_{2}\equiv 1$, that is $d^{n-1}(d+\omega(d'+dh'))e^{nh-\lambda_{2}z}\equiv p_{2}$. Obviously, $d^{n-1}(d+\omega(d'+dh'))\not\equiv0$. Otherwise, if $d^{n-1}(d+\omega(d'+dh'))\equiv0$,
we'll have
$p_{2}\equiv0$,  a contradiction.
So, by Lemma~\ref{lemma24} and Milloux's theorem, we have
\begin{eqnarray*}
  S\left(r,e^{h}\right)+nT\left(r,e^{h}\right)&=&T\left(r,e^{nh}\right)
  =T\left(r, \frac{p_{2}e^{\lambda_{2}z}}{d^{n-1}(d+\omega(d'+dh'))}\right)\\
  &\leq& T\left(r, p_{2}e^{\lambda_{2}z} \right)+T\left(r, d^{n-1}(d+\omega(d'+dh')) \right)+O(1)\\
  &=&O(r)+O(T(r,d)),
\end{eqnarray*}
which contradicts with the assumption that $\deg h=\sigma(f)>\max\{\sigma(d),1\}$.

If $f_{3}\equiv 1$, then by \eqref{theoremeqz12}, we have
$f_{1}+f_{2}\equiv 0$. That is
\begin{eqnarray*}
 d^{n-1}(d+\omega(d'+dh'))e^{nh}=p_{1}e^{\lambda_{1}z}.
\end{eqnarray*}
Following the similar reason as above, we can also get a contradiction.


If $\sigma(f)=1$, then we have $\sigma(d)<1=\deg h=\sigma(e^{(\lambda_{1}-\lambda_{2})z})$ and  $T(r)\geq T(r,f_{1}) = T\left(r,e^{(\lambda_{1}-\lambda_{2})z}\right)+S\left(r,e^{(\lambda_{1}-\lambda_{2})z}\right)$. Thus, by Milloux's theorem and Lemma~\ref{lemma23}, we get
\begin{eqnarray*}
 \frac{ N\left(r,\frac{1}{f_{2}} \right)}{T(r)}
  =\frac{O(T(r,d))+O(\log r)}{T\left(r,e^{(\lambda_{1}-\lambda_{2})z}\right)}
  \cdot\frac{T\left(r,e^{(\lambda_{1}-\lambda_{2})z}\right)}{T(r)}
  \to 0,
\end{eqnarray*}
and
\begin{eqnarray*}
 \frac{N\left(r,\frac{1}{f_{3}} \right)}{T(r)} \leq
 \frac{T(r,d)+S(r,d)+O(\log r)}{T\left(r,e^{(\lambda_{1}-\lambda_{2})z}\right)}\cdot\frac{T\left(r,e^{(\lambda_{1}-\lambda_{2})z}\right)}{T(r)}
  \to 0,
\end{eqnarray*}
as $ r\to \infty$.

Therefore, by using
Lemma~\ref{lemma27}, we can deduce that $f_{2}\equiv 1$ or $f_{3}\equiv 1$.

If $f_{2}\equiv 1$, that is
 \begin{eqnarray}\label{theoremeqz133}
 d^{n-1}(d+\omega(d'+dh'))e^{nh-\lambda_{2}z}= p_{2}.
\end{eqnarray}

We assert that $h'=\lambda_{2}/n$. Otherwise,
suppose that $h'\neq \lambda_{2}/n$, then from $\sigma(d')=\sigma(d)<1=\deg (nh-\lambda_{2}z)$, by using Lemma~\ref{lemma21} to \eqref{theoremeqz133}, we get $p_{2}\equiv0$, a contradiction.  Thus $h'=\lambda_{2}/n$. We set $h=\lambda_{2}z/n+B$, where $B$ is a
constant. By substituting it into \eqref{theoremeqz133}, we have
\begin{eqnarray}\label{theoremeqz13}
 d^{n-1}\left(d+\omega\left(d'+\frac{\lambda_{2}d}{n}\right)\right)
 =p_{2}e^{-nB}.
\end{eqnarray}
 Next, we assert that $d$ is a constant. Otherwise, if $d$ is a non-constant
 entire function, then from \eqref{theoremeqz13} we get that $0$ is a Picard exceptional value of $d$. Thus by Lemma~\ref{lemma28}, we have $d=e^{\alpha}$, where $\alpha$ is a non-constant polynomial, which contradicts with the assumption that $\sigma(d)<1$.
 So we have that $d$ is a non-zero constant, and \eqref{theoremeqz13} reduces to
 \begin{eqnarray*}
 d^{n}e^{nB}\left(1+\omega\frac{\lambda_{2}}{n}\right)
 =p_{2}.
\end{eqnarray*}
 Therefore,
 \begin{eqnarray*}
f=de^{h}=de^{B}e^{\lambda_{2}z/n}=\left(\frac{p_{2}n}{n+\omega\lambda_{2}} \right)^{\frac{1}{n}}e^{\frac{\lambda_{2}z}{n}}.
\end{eqnarray*}

Moreover, from $f_{2}\equiv 1$ and \eqref{theoremeqz12}, we also
 have $f_{1}+f_{3}\equiv0$. That is
\begin{eqnarray*}
p_{1}e^{\lambda_{1}z}
  = qd_{c}e^{Q+h_{c}},
\end{eqnarray*}
which implies that
\begin{eqnarray*}
 Q=\left(\lambda_{1}-\frac{\lambda_{2}}{n}\right)z+b_{1},
\end{eqnarray*}
where $b_{1}$ satisfies $p_{1}=q\left(\frac{p_{2}n}{n+\omega\lambda_{2}} \right)^{\frac{1}{n}}e^{\frac{\lambda_{2}c}{n}+b_{1}}$.

If $f_{3}\equiv 1$, then from \eqref{theoremeqz12} we have $f_{1}+f_{2}=0$. This gives that
\begin{eqnarray*}
 d^{n-1}(d+\omega(d'+dh'))e^{nh-\lambda_{1}z}=
 p_{1}.
\end{eqnarray*}
By using the similar methods as in the case $f_{2}\equiv 1$, we get
 \begin{eqnarray*}
f(z)=\left(\frac{p_{1}n}{n+\omega\lambda_{1}} \right)^{\frac{1}{n}}e^{\frac{\lambda_{1}z}{n}}.
\end{eqnarray*}
Furthermore, from $f_{3}\equiv 1$, we have
$qd_{c}e^{Q+h_{c}-\lambda_{2}z}\equiv p_{2}$. Thus we
can deduce
\begin{eqnarray*}
 Q=\left(\lambda_{2}-\frac{\lambda_{1}}{n}\right)z+b_{2},
\end{eqnarray*}
where $b_{2}$ satisfies $p_{2}=q\left(\frac{p_{1}n}{n+\omega\lambda_{1}} \right)^{\frac{1}{n}}e^{\frac{\lambda_{1}c}{n}+b_{2}}$. From the above discussion, the proof of Conclusion (2) is
complete.

\section{Proof of Theorem~\ref{thmABABnew3}.}

Suppose that $f$ is a transcendental entire solution with finite order of equation \eqref{thmABABnew2eq1addz2}.

If $\sigma(f)<1$, then following the
similar method as in the  proof of Case 1 of  Theorem~\ref{thmABABnew2}, we can get a contradiction.

If $\sigma(f)>1$. We denote $P=p_{1}e^{\lambda_{1} z}+p_{2}e^{\lambda_{2} z}$, and  rewrite
 \eqref{thmABABnew2eq1addz2}  as:
\begin{eqnarray}\label{theorem3eqz01}
 f^{3}+\omega f^{2}f'+(qf_{c})e^{Q}=P.
\end{eqnarray}

Differentiating \eqref{theorem3eqz01} yields
\begin{eqnarray}\label{theorem3eqz02}
 3f^{2}f'+\omega 2f(f')^{2}+ \omega f^{2}f''
 +Le^{Q}=P',
\end{eqnarray}
where $L=(qf_{c})'+Q'(qf_{c})$.

Eliminating $e^{Q}$ from \eqref{theorem3eqz01}  and \eqref{theorem3eqz02}, we have
\begin{eqnarray}\label{theorem3eqz033}
 fH=PL-P'(qf_{c}),
\end{eqnarray}
where
\begin{eqnarray*}
H=Lf^{2}+(\omega L-nqf_{c}) ff'-(n-1)\omega qf_{c}(f')^{2}-\omega qf_{c} f f''.
\end{eqnarray*}

If $H\equiv 0$, then from \eqref{theorem3eqz033} we have $PL-P'(qf_{c})\equiv0$. Similarly as in the proof of Subcase 2.1.2 of Theorem~\ref{thmABABnew2}, we get a contradiction. Therefore, $H\not\equiv 0$. From the assumptions that $f$ is entire and $q,Q$ are polynomials, we have that $H$ is an entire function. Since
both $PL-P'(qf_{c})$ and $H/f$ are all differential polynomials with meromorphic coefficients, similarly as in the proof of  Subcase 2.1.1 of Theorem~\ref{thmABABnew2},
by using Lemma~\ref{lemma25} to \eqref{theorem3eqz033}, we obtain that
\begin{eqnarray*}
T(r,H)=m(r,H)=S(r,f)+O(r).
\end{eqnarray*}
and
\begin{eqnarray*}
m(r,H/f)=S(r,f)+O(r).
\end{eqnarray*}

By observing we see that
the poles of $H/f$ all arise from the poles of $(f')^{2}/f$, thus
might arise from the zeros of $f$. Suppose that $z_{0}$ is a zero of $f$ with multiplicity $p$, then it
is a pole of $(f')^{2}/f$ with multiplicity $1$ when $p=1$; and
is a zero of $(f')^{2}/f$ with multiplicity $p-2$ when $p\geq2$.
So we have
\begin{eqnarray*}
T(r,H/f)&=&m(r,H/f)+N(r,H/f)=m(r,H/f)+N_{1)}(r,1/f)\\
&<&(\kappa+o(1))T(r,f)+S(r,f)+O(r).
\end{eqnarray*}

Therefore, by $H\not\equiv 0$ we have
\begin{eqnarray*}
T(r,f)&=&T(r,1/f)+O(1)\leq T\left(r,\frac{H}{f}\right)+T\left(r, \frac{1}{H}\right)+O(1)\\
&<&(\kappa+o(1))T(r,f)+S(r,f)+O(r),  \quad  0\leq\kappa<1.
\end{eqnarray*}
So we have
\begin{eqnarray*}
T(r,f)=S(r,f)+O(r),
\end{eqnarray*}
which contradicts with the assumption that $\sigma(f)>1$.

If $\sigma(f)=1$, then following the
similar method as in the  proof of Case 3 of  Theorem~\ref{thmABABnew2}, we can get that $\sigma(f)=\deg Q=1$.


\noindent \textbf{Acknowledgements.} Thanks are due to the referee for valuable comments and suggestions.

\def\cprime{$'$}
\providecommand{\bysame}{\leavevmode\hbox to3em{\hrulefill}\thinspace}
\providecommand{\MR}{\relax\ifhmode\unskip\space\fi MR }
\providecommand{\MRhref}[2]{%
  \href{http://www.ams.org/mathscinet-getitem?mr=#1}{#2}
}
\providecommand{\href}[2]{#2}

\end{document}